\documentclass[11pt]{amsart} 
\usepackage{amsmath,amssymb}
\numberwithin{equation}{section} 
\usepackage{xcolor}

\newtheorem{theorem}{Theorem}[section]

\newtheorem*{theoremA}{Theorem A}
\newtheorem*{theoremB}{Theorem B}

\newtheorem*{theoremC}{Theorem C}

\newtheorem{thm}[theorem]{Theorem}
\newtheorem{lemma}[theorem]{Lemma}
\newtheorem{lem}[theorem]{Lemma}
\newtheorem{prop}[theorem]{Proposition}
\newtheorem{proposition}[theorem]{Proposition}

\newtheorem{conj}[theorem]{Conjecture}
\newtheorem{pbm}[theorem]{Problem}

\theoremstyle{definition}

\newtheorem{defi}[theorem]{Definition}
\newtheorem{example}[theorem]{Example}
\newtheorem{exa}[theorem]{Example}

\newtheorem{remark}[theorem]{Remark}
 \newtheorem*{ackn}{Acknowledgements}

\def \R{\mathbb{R}}
\def \C{\mathbb{C}}
\def \B{\mathbb{B}}

\def \N{\mathbb{N}}

  \newcommand{\f}{\varphi}
 
 \newcommand{\p}{\psi}

  \newcommand \e {\varepsilon}
 \newcommand \la {\lambda}

 \usepackage{hyperref}
\hypersetup{
    unicode=false,        
    pdftoolbar=true,      
    pdfmenubar=true,       
    pdffitwindow=false,     
    pdfstartview={FitH},    
    pdftitle={},    
    pdfauthor={},     
    colorlinks=true,       
   linkcolor=black,          
    citecolor=black,        
    filecolor=black,      
    urlcolor=black}

\begin{document}

\title[High Energy Classes]{High Energy plurisubharmonic classes}

\date{\today}

\author{Vincent Guedj}
\address{Institut Universitaire de France \& Institut de Mathématiques de Toulouse\\
  Université de Toulouse, CNRS\\
  118 route de Narbonne, F-31400 Toulouse\\
  France}
  \email{\href{vincent.guedj@math.univ-toulouse.fr}{vincent.guedj@math.univ-toulouse.fr}}
\urladdr{\href{https://www.math.univ-toulouse.fr/~guedj/}{https://www.math.univ-toulouse.fr/~guedj/}}

 \author{Ahmed Zeriahi}
 \address{Institut de Mathématiques de Toulouse\\
  Université de Toulouse, CNRS\\
  118 route de Narbonne, F-31400 Toulouse\\
  France}
  \email{\href{ahmed.zeriahi@math.univ-toulouse.fr}{ahmed.zeriahi@math.univ-toulouse.fr}}

\begin{abstract} 
Let $\Omega \Subset \C^n$ be  a bounded strongly pseudoconvex domain. 
For any concave increasing weight $\chi : \R^- \longrightarrow \R^-$ such that $\chi(0) = 0$, we 
introduce and study finite energy classes $\mathcal E_\chi(\Omega)$ of plurisubharmonic functions,
using the Orlicz space formalism.
We investigate the range of the Monge-Ampère operator on these classes, and conjecture that
this should lead to an integral characterization of the image of bounded 
plurisubharmonic functions, an open problem since the birth of pluripotential theory more than forty years ago.
\end{abstract}

\maketitle


\section*{Introduction}

The classical Dirichlet energy is a fundamental tool in the study of the Laplace equation. 
A weighted and non linear version  has been introduced by Cegrell in \cite{Ceg98}, in order to study the  complex Monge-Amp\`ere operator.
Extending these ideas to the context of compact K\"ahler manifolds \cite{GZ07}
has revealed extremely useful; we refer the interested reader to 
\cite{GZbook,Bou18} for some overview.

In the last decade most works have been devoted to fine properties of finite energy
classes ${\mathcal E}_{\chi}$ with respect to weights $\chi$ that have slow growth at infinity
(see e.g. \cite{Dar24}), as the latter 
allow one to characterize probability measures that do not charge pluripolar sets
(see \cite[Theorem A]{GZ07}).
In this article we take a closer look at  weights with fast growth at infinity,
studying the corresponding {\it high energy classes}.
We restrict to the case of euclidean domains for simplicity.

 \smallskip

Let $\Omega=\{ z \in \C^n, \, \rho(z)<0 \}$ be a bounded strongly pseudoconvex domain,
 with $\rho:\C^n \rightarrow \R$ a smooth strictly plurisubharmonic function.
 Let ${\mathcal T}(\Omega)$ denote the set of {\it test functions}: these are bounded plurisubharmonic functions
 $u$ in $\Omega$ such that $u_{|\partial \Omega}=0$ and $\int_{\Omega} MA(u)<\infty$,
 where $MA(u)$ denotes the complex Monge-Amp\`ere measure of $u$ (in the sense of Bedford-Taylor).
  Let $h: \R^+ \rightarrow \R^+$ be  a convex increasing function such that 
 $h(0)=0$, $h'(0)=1$ and set $\chi(t)=-h(-t)$;
 following \cite{GZ07}  we let ${\mathcal W}^+$ denote the set of all such weights $\chi$.
   The finite energy class ${\mathcal E}_{\chi}(\Omega)$ is the set of decreasing limits $\f=\lim \searrow u_j$
 of test functions $u_j$ such that  $\sup_j E_{\chi}(u_j)<+\infty$, where
 $$
 E_{\chi}(u_j):=\inf \left\{ \la>0, \; \int_{\Omega} (-\chi) (u_j/\la) MA(u_j) \leq 1 \right \}.
 $$
 
 We establish in this article some key properties of the classes ${\mathcal E}_{\chi}(\Omega)$.

\begin{theoremA}
Fix a weight $\chi \in {\mathcal W}^+$.
\begin{enumerate}
\item The complex Monge-Amp\`ere operator $MA$ is well defined on ${\mathcal E}_{\chi}(\Omega)$ and
$$
{\mathcal E}_{\chi}(\Omega)=\left\{ \f \in {\mathcal E}^1(\Omega); \int_{\Omega} (-\chi) (\f/\la) MA(\f) <+\infty
\text{ for some } \la>0 \right\}.
$$

\item The energy $E_{\chi}$ is semi-continuous along decreasing sequences in ${\mathcal E}_{\chi}(\Omega)$.

\smallskip

\item If $\f, \p$ belong to ${\mathcal E}_{\chi}(\Omega)$ with $\f \leq \p$, then  $0 \leq E_{\chi}(\p) \leq 2^{n+1} E_{\chi}(\f)$.

\smallskip

\item A function $\f$ belongs to ${\mathcal E}_{\chi}(\Omega)$ if and only if there exists $\la,M>0$ and a decreasing limit 
$\lim \searrow u_j=\f$ of test functions $u_j$ such that 
 $$
 \int_0^{+\infty} t^n \chi'(-t/\la) {\rm Cap}(u_j<-t) dt \leq M,
 $$
 where ${\rm Cap}$ denotes Bedford-Taylor's Monge-Amp\`ere capacity.
 
 \item The class ${\mathcal E}_{\chi}(\Omega)$ is a convex cone with strong $L^1$ compactness properties.
\end{enumerate}
\end{theoremA} 

This result is a combination of Theorem \ref{thm:energy-approximation}, 
Lemma \ref{lem:fdtl}, Propositions \ref{pro:highenergyequality}, \ref{prop:cap-characterization}
and \ref{pro:energycpct}.
A different notion of high energy classes had been previously introduced in \cite{GZ07,BGZ09}, but it has poorer properties and is less suited to describe the range of the Monge-Amp\`ere operator.
An important new idea of this article lies indeed in the use of the Orlicz formalism to study the latter.

We also establish a subextension property (Proposition \ref{pro:subext}) and a Moser-Trudinger inequality (Theorem \ref{thm:mosert}),
which yield strong integrability properties and generalize  results
of Cegrell-Zeriahi \cite{CZ03} and Di Nezza-Guedj-Lu \cite{DGL21}.

\smallskip

We then study the range of the complex Monge-Amp\`ere operator on the classes ${\mathcal E}_{\chi}(\Omega)$, 
trying to generalize the work of Cegrell \cite{Ceg98} who dealt with weights
$\chi(t)=-\frac{(-t)^p}{p}$ with polynomial growth.
We provide a new proof of Cegrell's main results and propose a conjectural description
of the range $MA({\mathcal E}_{\chi}(\Omega))$ for weights with arbitrary growth.
Setting
$$
\kappa(\mu,\chi):=\sup_{\f \in {\mathcal E}_{\chi}(\Omega) \setminus \{0\}} \frac{||\f||_{L_{\chi}(\mu/2)}}{ ||\f||_{L_{\chi}(\mu)}},
$$
when ${\mathcal E}_{\chi}(\Omega) \subset L_{\chi}(\mu)$, we conjecture that this quotient is always less than $1$.

 \begin{theoremB}
Fix $\mu$ a positive Radon measure and $\chi \in {\mathcal W}^+$.
\begin{enumerate}
\item If there exists $a \in (0,1)$ and $C>0$ such that for all $\f \in {\mathcal E}_{\chi}(\Omega)$,
\begin{equation}
\tag{$I(a,\chi)$}
||\f||_{L_{\chi}(\mu)} \leq a E_{\chi}(\f)+C, 
\end{equation}
then there exists a unique $\psi \in \mathcal E_\chi(\Omega)$ such that $\mu = MA(\psi)$.

\item If there exists   $\psi \in \mathcal E_\chi(\Omega)$ such that $\mu = MA(\psi)$,
then for any $a>1$ there exists $C>0$ such that $I(a,\chi)$ holds for all $\f \in {\mathcal E}_{\chi}(\Omega)$.

\item If $\kappa(\mu,\chi)<1$ and there exists   $\psi \in \mathcal E_\chi(\Omega)$ such that $\mu = MA(\psi)$
then there exists $a \in (0,1)$ and $C>0$ such that $I(a,\chi)$ holds for all $\f \in {\mathcal E}_{\chi}(\Omega)$.
\end{enumerate}
\end{theoremB} 

This statement is a combination of Proposition \ref{pro:quantitative}, Theorems \ref{thm:quantitative} and \ref{thm:range}.
We finally apply our analysis to 
a long standing open problem: does there exist an integral 
characterization of the range of the complex Monge-Amp\`ere operator acting on {\it bounded  plurisubharmonic functions} ?

 \begin{theoremC}
Assume that $\kappa(\mu,\chi)<1$
for all weights $\chi \in {\mathcal W}^+$ and positive Radon measures $\mu$ with ${\mathcal E}_{\chi}(\Omega) \subset L_{\chi}(\mu)$.
Then the following are equivalent:
\begin{enumerate}
\item $\mu$ is the Monge-Amp\`ere measure of a bounded plurisubharmonic function
with zero boundary values;
\item for any $\chi \in {\mathcal W}^+$, there exists $a \in (0,1)$ and $C>0$ such that
$I(a,\chi)$ holds.
\end{enumerate}
\end{theoremC} 
 
Polynomial weights provide interesting examples for which $\kappa(\mu,\chi)<1$.
The case of general weights is  a problem of independent interest.
We refer the interested reader to \cite{Bed87} and \cite[Questions 13,14,15,17]{DGZ16} for more motivation.

\begin{ackn} 
We thank S.Dinew for useful remarks.
VG is supported by the Institut Universitaire de France
and the fondation Charles Defforey. 
\end{ackn}

\section{Preliminaries} 

In the whole article we fix 
$\Omega=\{ z \in \C^n, \, \rho(z)<0 \}$
a bounded strongly pseudoconvex domain,
 where $\rho:\C^n \rightarrow \R$ is a smooth strictly plurisubharmonic function.
Many properties are valid when $\Omega$ is rather a bounded {hyperconvex} domain
 of a Stein space,
but we stick to this   setting to simplify the discussion.

\subsection{Orlicz  Spaces} 

\subsubsection{The space of weights ${\mathcal W}^+$}

\begin{defi}
 A Young function is a convex increasing function
  $h : \R^+ \longrightarrow \R^+$  such that $h(0) =0$
 and  $\lim_{s \to + \infty} \frac{h(s)}{s} = + \infty$.
\end{defi}

Fix $\mu$ a positive Borel measure 
in $\Omega$, and  $h$ a Young function.

\begin{defi}
The Orlicz space $ L_h(\mu)=L_h(\Omega,\mu)$  is the set of   measurable  functions  
$f : \Omega \longrightarrow \bar \R $ such that there is $\lambda > 0$ with
$
  \int_{\Omega} h\left( \frac{|f|}{\lambda} \right) d \mu < + \infty.
$
\end{defi}

This is a vector space endowed with a natural {\it Luxembourg norm},  defined by
$$
\Vert f\Vert_{L_h(\mu)}  := \inf \left\{ \lambda > 0 ; \int_\Omega h(\vert f\vert\slash \lambda) d \mu \leq 1 \right\}\cdot
$$
The classical Lebesgue spaces $L^p(\mu)$ correspond to 
$h(t)=t^p/p$.
Observe that 
\begin{itemize}
\item the normalization $h(0)=0$ ensures that $\Vert f\Vert_{L_h(\mu)} = 0$ if and only if  $f = 0$;
\item  if $\Vert f\Vert_{L_h(\mu)} > 0$ then 
$
\int_\Omega h(\vert f\vert\slash \Vert f\Vert_{h,\mu}) d \mu = 1.
$
\end{itemize}

\begin{defi}
Fix $h$ a Young function. We consider 
$\chi : \R^- \longrightarrow \R^-$ defined by $\chi(t)=-h(-t)$.
We let  ${\mathcal W}^+$ denote the set of all  weights $\chi$
such that $\chi'(0)=1$.
\end{defi}

A weight $\chi \in {\mathcal W}^+$ is a 
concave increasing  function such that $\chi(0) = 0$
 and $\lim_{t \to - \infty} \frac{\chi (t)}{t} = + \infty$.
The normalization  $\chi'(0)=1$ ensures that $\chi(t) \leq t$.

\subsubsection{Extra growth}

We shall need the following classical lemma.

\begin{lem} \label{lem:unpeu+}
Let $f:\Omega \rightarrow \R^+$ be a Lebesgue measurable function.
Let $h$ be a Young function such that $\int_{\Omega} h \circ f d\mu<+\infty$. 
For all $\kappa>1$, there exists  a Young function $\hat{h}$ such that $\hat{h}(t)/h(t) \rightarrow +\infty$ and 
$$
\int_{\Omega} \hat{h} \circ f d\mu \leq \kappa \int_{\Omega} h \circ f d\mu .
$$
\end{lem}

\begin{proof}
It suffices to prove the following related result: if $f \in L^1(\mu)$ then there exists 
a Young function $h$ such that $\int h \circ f d \mu \leq \kappa \int f d\mu$.
By Fubini theorem
$$
\int h \circ f d \mu=\int_0^{+\infty} h'(t) F(t) dt,
\; \text{ with } \; F(t)=\mu(f>t).
$$
Given $F: \R^+\rightarrow \R^+$ a decreasing function such that $\int_0^{+\infty} F(t) dt<+\infty$, the problem thus boils down
to finding an increasing function $a: \R^+\rightarrow \R^+$ such that $a(t) \rightarrow +\infty$ as $t \rightarrow +\infty$,
and $\int_0^{+\infty} a(t) F(t) dt \leq \kappa \int_0^{+\infty} F(t) dt$.

Fix $\kappa>1$, $0<q<1$ such that $\kappa(1-q)=1$ and
set $I=\int_0^{+\infty} F(t) dt$. 
Since $\int_N^{+\infty} F(t) dt \rightarrow 0$ as $N \rightarrow +\infty$, we can pick $N_k$ an increasing sequence
of integers such that $N_0=0$, $N_{k+1} \geq 2 N_k$, and
$
\int_{N_k}^{+\infty} F(t) dt  \leq q^{2k} I.
$
We set $a(t)=q^{-k}$ for $N_k \leq t<N_{k+1}$. Then 
$
 \int_{N_k}^{N_{k+1}} a(t) F(t) dt  \leq q^{k} I,
 $
 hence
 $
 \int_{0}^{+\infty} a(t)F(t) dt  \leq \sum_{k=0}^{+\infty} q^{k} I= \kappa \int_0^{+\infty} F(t) dt.
 $
The conclusion follows, setting $h(t)=\int_0^t a(s)ds$.
\end{proof}


\subsection{Range of the Monge-Amp\`ere operator} \label{sec:prelim2}

\subsubsection{Bedford-Taylor theory}

We let $PSH(\Omega)$ denote the set of functions that are defined and plurisubharmonic (psh) in $\Omega$;
functions in $PSH(\Omega)$ are locally integrable
with respect to the euclidean Lebesgue measure $dV$.
We endow  $PSH(\Omega)$ with the $L^1_{loc}$-topology.
Skoda's  integrability result \cite[Theorem 2.50]{GZbook} asserts that psh functions
are actually {\it exponentially integrable}, namely 
$$
PSH(\Omega) \subset L_h(K,dV)
$$
for each compact subset $K \subset \Omega$, with $h(t)=e^t-1$.

\smallskip

When $u \in PSH(\Omega)$ is smooth, its complex Monge-Amp\`ere measure
$$
MA(u)=(dd^c u)^n=c_n \det ({\rm Hess}(u)) dV
$$
 is a smooth positive measure.
 Here $d=\partial+\overline{\partial}$ and 
$d^c=\frac{1}{2i\pi}(\partial-\overline{\partial})$ are real operators, 
$dd^c =\frac{i}{\pi}\partial\overline{\partial}$ and $c_n>0$.
When $u \in PSH(\Omega)$ is merely {\it locally bounded}, its Monge-Amp\`ere measure
$MA(u)=(dd^c u)^n$ is a well defined positive Radon measure in $\Omega$,
as shown by Bedford and Taylor in their foundational work \cite{BT82}.

\begin{defi}
${\mathcal T}(\Omega)$ denotes the 
set of bounded negative plurisubharmonic
functions $u$ in $\Omega$ such that
 $\lim _{z\to \zeta} u (z) = 0,$ 
$\forall \zeta \in \partial \Omega,$ and $\int_\Omega(dd^c u )^n <+\infty$.
\end{defi}

Functions in this convex cone play the role of {\it test functions} for the Monge-Amp\`ere operator.
It is useful  to allow for tests that are not necessarily smooth.
While $MA(u)$ has finite mass in $\Omega$ when $u \in {\mathcal C}^{\infty}(\overline{\Omega})$, this is not necessarily the case
when $u$ is merely bounded, as one can easily check on radial examples.

\begin{pbm} \label{pbm:bedford}
\cite[p.70]{Bed87}
Can one characterize the range of the complex Monge-Amp\`ere operator $MA$ acting
on $PSH(\Omega) \cap L^{\infty}_{loc}$ ?
\end{pbm}

Forty years after the seminal work \cite{BT82} this problem remains open
(see \cite[Question 13]{DGZ16}) and is the main motivation
for the present work.

\subsubsection{Dirichlet problem}

Solving the Dirichlet problem is a  fundamental question in PDEs.
Given $\mu$ a positive Radon measure, one seeks  a function $u \in PSH(\Omega)$
such that $u_{|\partial \Omega}=0$ and $MA(u)=\mu$.
When $\mu=f dV$ is absolutely continuous with respect to the Lebesgue measure,
Kolodziej has provided in \cite{Kol05}  criteria on the density $f$
ensuring that there exists a unique bounded solution.
Studying this problem for mildly unbounded functions is   very useful,
as shown by Cegrell.

\subsubsection*{Cegrell classes}

Cegrell \cite{Ceg98, Ceg04} has introduced   several
classes of  unbounded plurisubharmonic functions on which 
the   operator $MA$ is well defined: 

\begin{itemize}
\item  ${\rm DMA}(\Omega)$ is the set of psh functions
 $u$ such that for all $z_0 \in \Omega$, there exists a neighborhood $V_{z_0}$ of
$z_0$ and $u_j \in {\mathcal T}(\Omega)$ a decreasing sequence which
converges to  $u$ in $V_{z_0}$ and satisfies
$\sup_j \int_{\Omega} (dd^c u_j)^n <+\infty$.

\item 
$(dd^c \cdot )^n$ is well defined
on $DMA(\Omega)$ and continuous under decreasing limits;
the class $DMA(\Omega)$ is stable under taking maximum
(see also \cite{Blo06}).


\item The class ${\mathcal E}^p(\Omega)$ 
is the set of psh functions $u$ for which there exists a sequence of
functions $u_j \in {\mathcal T}(\Omega)$ decreasing towards $u$ in all of $\Omega$, and
so that $\sup_j \int_{\Omega} (-u_j)^p (dd^c u_j)^n<+\infty$.
The weighted energy of $u$ by 
\[
E_p(u):= \int_{\Omega} (-u)^p (dd^c u)^n<+\infty.
\]  
\item  If $u\in \mathcal{E}^p(\Omega)$ for some $p>0$ then $(dd^c u)^n$ vanishes on all pluripolar sets \cite[Theorem 2.1]{BGZ09}. 
 Also, note that 
$
{\mathcal T}(\Omega) \subset {\mathcal E}^p(\Omega)\subset {\rm DMA}(\Omega). 
 $
\end{itemize}

 A beautiful result of Cegrell characterizes the range of the operator $MA$:
 
 \begin{thm}\cite[Theorem 6.2]{Ceg98} \label{thm:ceg}
 Let $\mu$ be a positive Radon measure in $\Omega$. The following are equivalent:
 \begin{enumerate}
 \item there is a unique function $u \in {\mathcal E}^p(\Omega)$ such that   $MA(u)=\mu$;
 \item ${\mathcal E}^p(\Omega) \subset L^p(\Omega,\mu)$.
 \end{enumerate}
  \end{thm}
 
 Understanding  the range of the complex Monge-Amp\`ere operator 
on bounded plurisubharmonic functions requires one to extend this study to the
case of   finite energies with respect to a weight $\chi$ having non polynomial growth.

 \subsubsection{Maximum principles}

\begin{theorem} \cite[Theorem 4.1]{BT82} \label{thm:comparisonprinciple}
Fix $u,v \in {PSH}(\Omega) \cap L^{\infty} (\Omega)$ such that $\liminf_{z \to \zeta} (u - v) (z) \geq 0$. Then  
$$
\int_{\{u<v\}} MA(v) \leq \int_{\{u<v\}} MA(u).
$$
\end{theorem}

\noindent This  {\it comparison principle} is an integral form of the classical maximum principle.
It has played  a central role  in the applications
of pluripotential theory to K\"ahler geometry (see \cite{GZbook}).
The situation is more delicate in ${\mathcal E}^1(\Omega)$, but 
the following  weaker version  will be sufficient for us here (see \cite[Corollary 6]{Xing08}). 

\begin{theorem} \label{thm:comparisonprinciple2}
Fix $u,v \in {\mathcal E}^1(\Omega)$. Then  
\begin{enumerate}
\item One has ${\bf 1}_{\{u>v\}} MA(\sup\{u,v\}) = {\bf 1}_{\{u>v\}} MA(u)$.

\smallskip

\item If $MA(u) \leq MA(v)$ then $u \geq v$.
\end{enumerate}
\end{theorem}

\subsection{Monge-Amp\`ere capacity}

We let  ${\rm Cap}(\cdot) := \text{Cap}(\cdot,\Omega)$ denote the Monge-Amp\`ere capacity 
introduced in \cite{BT82}: 

\begin{defi} \label{def:cap}
If $E \subset \Omega$ is a Borel subset, then
\[
{\rm Cap}(E,\Omega) := \sup \left \{ \int_E (dd^c u)^n \; | \; u \in PSH(\Omega), \ -1\leq u \leq 0 \right \}. 
\]
\end{defi}

It follows from the work of Bedford-Taylor \cite{BT82} that ${\rm Cap}$ is a {\it Choquet capacity}, i.e. 
it is a set function satisfying the following properties:

\smallskip

$(i)$ ${\rm Cap}(\emptyset) = 0$;

\smallskip

$(ii)$ ${\rm Cap}$ is monotone, i.e.
$
 A \subset B \subset \Omega \Longrightarrow 0 \leq {\rm Cap}(A) \leq {\rm Cap}(B);
$

$(iii)$  if $(A_n )_{n \in \N}$ is a non-decreasing sequence of subsets of $\Omega$, then
$$
{\rm Cap} (\cup_n A_n) = \lim_{n \to + \infty} {\rm Cap}(A_n ) = \sup_n {\rm Cap}(A_n );
$$

$(iv)$ if $(K_n )$ is a non-increasing sequence of compact subsets of
$\Omega$,
$$
{\rm Cap}(\cap_n K_n )  = lim_{n \to + \infty} {\rm Cap}(K_n) = \inf_n {\rm Cap}(K_n).
$$

$(v$)  
if $(A_n)_{n \in \N}$ is any sequence of subsets of $\Omega$, then
$
{\rm Cap}(\cup_n A_n) \leq \sum_n {\rm Cap} (A_n).
$

\smallskip

 In particular  
$
{\rm Cap} (B) = \sup \{ {\rm Cap} (K) ; K \,  \text{compact} \,  \subset B\}
$
for every Borel set $B$.
  This is a special case of Choquet's capacitability theorem \cite{Cho54}.

\medskip

We  need the following result  (\cite[Prop. 2]{Kol96} and \cite[Prop. 3.4]{CKZ05}).

\begin{lemma}  \label{lem:classic}
Let  $\varphi \in \mathcal T (\Omega)$. Then for any $s, t > 0$ we have
\begin{equation*}
t^{n } \mathrm{Cap} (\{\varphi < - s - t\}) \leq  \int_{\{\varphi <  - s\}}   (dd^c \varphi)^n 
\leq s^{n} \mathrm{Cap}  (\{\varphi <  - s\}).
\end{equation*}
\end{lemma}

\section{High energy classes}

In the whole section we fix a weight  $\chi \in {\mathcal W}^+$.

\subsection{Basic properties}

\begin{defi}
The weighted $\chi$-energy  of a function $\varphi \in \mathcal T (\Omega)$ is 
$$
E_\chi(\varphi) :=  
\inf \left\{ \lambda > 0 ; \int_\Omega h(\vert \varphi \vert\slash \lambda) d \mu_\varphi \leq 1 \right\},
$$
 where $ \mu_\varphi = (dd^c \varphi)^n$ denotes the Monge-Ampère measure of $\varphi$. 
\end{defi}

Although this definition is inspired by the one for Orlicz spaces, let us stress that $E_\chi(\varphi)$ is not a norm
on $\mathcal T (\Omega)$.

\begin{defi}
The class $\mathcal E_\chi(\Omega)$ is the set of $\varphi \in PSH(\Omega)$ such that there exists  a decreasing sequence $(\varphi_j)$ in $\mathcal T(\Omega)$ converging to $\varphi$ with
$\sup_j E_\chi(\varphi_j) < + \infty$.
\end{defi}

It follows from the definition that $\mathcal E_\chi(\Omega) \subset \mathcal E_{\tilde{\chi}}(\Omega)$ 
if $\chi \leq \tilde{\chi}$.
In particular $\mathcal E_\chi(\Omega) \subset \mathcal E^1(\Omega)$ since $\chi'(0)=1$, hence
the Monge-Amp\`ere measure $ \mu_\varphi := (dd^c \varphi)^n$
is well defined and we can consider the following energy functional:

\begin{defi}
For $\varphi \in \mathcal E^1 (\Omega)$ we set
$
E_\chi(\varphi) := \inf \left\{ \lambda > 0 ; \int_\Omega h(-\varphi \slash \lambda) d \mu_\varphi \leq 1 \right\}.
$
\end{defi}

We are going to show that $\mathcal E_\chi(\Omega)$ is precisely the subset 
of those functions in ${\mathcal E}^1(\Omega)$  that have finite $E_{\chi}$-energy.

\begin{lemma}  \label{lem:fdtl0}
Fix $\varphi \in {\mathcal E}^1(\Omega)$ and $\alpha >0$.  Then 
  $$
  E_\chi(\alpha \varphi) \leq \max \{\alpha, \alpha^{n+1}\} E_\chi(\varphi).
  $$
\end{lemma}

\begin{proof} 
Set $\lambda = \alpha E_\chi(\varphi)$. Then 
$$
\alpha^{-n}  \int_\Omega h(-\alpha \varphi\slash \lambda) MA(\alpha \varphi) =  \int_\Omega h(-\alpha \varphi \slash \lambda) MA(\varphi) = 1.
$$

If $\alpha \leq 1$ then $\int_\Omega h(-\alpha \varphi\slash \lambda) MA(\alpha \varphi)  = \alpha^{n} \leq 1$, 
hence we obtain $\lambda \geq E_\chi(\alpha \varphi)$ and $ \alpha E_\chi(\varphi) \geq  E_\chi(\alpha \varphi)$.
On the other hand if $\alpha > 1$, the convexity of $h$ yields
$$
\int_\Omega h(- \alpha \varphi\slash \lambda \alpha^{n}) MA(\alpha \varphi) \leq \alpha^{-n} \int_\Omega h(- \alpha \varphi\slash \lambda) MA(\alpha \varphi) \leq 1.
$$
This implies that $ \lambda \alpha^{n} \geq E_\chi(\alpha \varphi)$, hence $\alpha^{n+1} E_\chi(\varphi) \geq
E_\chi(\alpha \varphi)$. 
\end{proof}

\subsection{Continuity of the weighted energy}

\begin{theorem}  \label{thm:energy-approximation}
Fix   $\varphi \in  \mathcal E_\chi(\Omega)$.  
\begin{enumerate}
\item If  a sequence $(\varphi_j)_{j \in \N} \subset \mathcal E_\chi(\Omega)$ 
decreases to $\varphi$, then  
$$
E_\chi(\varphi) \leq \liminf_{j \to + \infty} E_\chi(\varphi_j).
$$

\smallskip

\item There exists  $(\widetilde{\varphi}_j)_{j \in \N} \subset \mathcal T(\Omega)$ decreasing to $\varphi$
such that  the measures $(dd^c \widetilde{\varphi}_j)^n$ converge to
$(dd^c \f)^n$ in the strong sense of Borel measures, and
$$
E_\chi(\varphi) = \lim_{j \to + \infty} E_\chi(\widetilde{\varphi}_j).
$$
\end{enumerate} 
\end{theorem}

\begin{proof}  
We can assume without loss of generality  $\liminf_{j \to + \infty} E_\chi(\varphi_j) < + \infty$.
Fix $\lambda > \liminf_{j \to + \infty} E_\chi (\varphi_j)$ and  $(j(k))_{k \in \N}$
 such that $ E_\chi (\varphi_{j(k)})  < \lambda$ for all $k$. Thus
$$
\int_\Omega h(-\varphi_{j(k)} \slash  \lambda)MA(\varphi_{j(k)})  
\leq \int_\Omega h(-\varphi_{j(k)} \slash  E_\chi (\varphi_{j(k)})) MA(\varphi_{j(k)})\leq 1.
$$
Since $k \longmapsto h(-\varphi_{j(k)} \slash  \lambda)$ is an increasing sequence of lower semi-continuous functions in $\Omega$ converging to $h(-\varphi \slash  \lambda) $ and $MA(\varphi_{j(k)}) \to MA(\varphi)$ 
 in the sense of Radon measures as ${k \to +\infty}$, we obtain
$$
\int_\Omega h(-\varphi \slash  \lambda) MA(\varphi) \leq \liminf_{k \to + \infty} \int_\Omega h(-\varphi_{j(k)}\slash  \lambda)MA(\varphi_{j(k)})  \leq 1.
$$
Thus $\lambda \geq E_\chi(\varphi)$, showing that
$
E_\chi (\varphi) \leq \liminf_{j \to + \infty} E_\chi (\varphi_j).
$

\medskip

We now prove 2), adapting some arguments from \cite[Theorem 5.6]{Ceg98}.
For $j \in \N$ we consider  $A_j := \{z \in \Omega ; \rho (z) <  - 2^{-j}, \varphi > - j\}\cdot$
This is a relatively compact Borel set in $\Omega$. 
It follows from Theorem \ref{thm:comparisonprinciple2} that 
$$
{\bf 1}_{A_j} (dd^c \varphi)^n = {\bf 1}_{A_j} (dd^c \max \{\varphi,-j\})^n
$$
Now \cite[Theorem 6.1]{Ceg98} ensures the existence of 
 $\widetilde{\varphi}_j \in \mathcal E^1(\Omega)$ such that
$$
(dd^c \widetilde{\varphi}_j)^n = {\bf 1}_{A_j} (dd^c \varphi)^n.
$$
Observe that $\widetilde{\varphi}_j \geq \max \{\varphi,-j\}$ (by Theorem \ref{thm:comparisonprinciple2}), 
hence $\widetilde{\varphi}_j$ is bounded. 

We claim that $\widetilde{\varphi}_j \in {\mathcal T}(\Omega)$.
Indeed set $M_j= j 2^{j+1}$, $K_j= \{z \in \Omega ; \rho (z) \leq - 2^{-j} \}$
 and $ \Omega_j := \{\rho<-2^{-j-1}\} \supset K_j$. 
Since $\widetilde{\varphi}_j \geq -j $ in $\Omega$, 
we obtain 
$$
M_j \rho \leq  -  M_j 2^{-j+1} \leq - j \leq  \widetilde{\varphi}_j
\; \; \text{ in } \; \; \partial \Omega_j.
$$
Set   
$\psi_j := \sup \{\widetilde{\varphi}_j, M_j \rho\}$ in $\Omega \setminus \Omega_j$ 
and $\psi_j := \widetilde{\varphi}_j$ in $ \Omega_j  \supset K_j$.
By the gluing principle, $\psi_j \in \mathcal T(\Omega)$.
Since  $\psi_j = \widetilde{\varphi}_j $ is a neighborhood of $\Omega_j$, we obtain 
$(dd^c \psi_j)^n \geq (dd^c \widetilde{\varphi}_j)^n$ in $\Omega$,
hence $\psi_j \leq \widetilde{\varphi}_j$ in $\Omega$
by Theorem \ref{thm:comparisonprinciple2}.

It follows from Theorem \ref{thm:comparisonprinciple2} again that   $ j \mapsto \widetilde{\varphi}_j$ is decreasing 
with  $\widetilde{\varphi}_j  \geq \varphi$ in $\Omega$. 
Setting $\lambda = E_\chi(\varphi)$ we obtain
$$
\int_\Omega -\chi(\widetilde{\varphi}_j \slash \lambda) (dd^c \widetilde{\varphi}_j )^n \leq \int_\Omega -\chi(\varphi \slash \lambda) (dd^c \varphi)^n =1.
$$ 
Hence $E_\chi(\widetilde{\varphi}_j ) \leq \lambda = E_\chi(\varphi)$ and 
$
\limsup_{j \to + \infty} E_\chi(\widetilde{\varphi}_j ) \leq E_\chi(\varphi).
$

On the other hand, the sequence $(\widetilde{\varphi}_j )_{j \in \N}$ decreases to a plurisubharmonic function $\psi $ such that $\psi \geq \varphi$ in $\Omega$. The  uniform bound on the energies shows that $\psi \in   \mathcal E_\chi (\Omega)$.
The continuity of $MA$ along decreasing sequences yields
$$
(dd^c \psi)^n = \lim_{j \to + \infty} (dd^c \widetilde{\varphi}_j )^n =  (dd^c \varphi)^n
$$
We conclude that $\psi = \varphi$  
by uniqueness of solutions (Theorem \ref{thm:comparisonprinciple2}).
Thus  $\widetilde{\varphi}_j$ converges to $\varphi$.
The proof is complete as $E_\chi(\varphi) \leq \liminf_{j \to + \infty} E_\chi(\widetilde{\varphi}_j )$ by 1).
\end{proof}

\begin{remark} \label{rem:keylemE1}
Using the canonical approximants constructed above, one can show that Lemma \ref{lem:classic} is valid for 
a function $\f$ that merely belongs to  ${\mathcal E}^1(\Omega)$.
\end{remark}

We now establish the appropriate high-energy version of the ``fundamental inequality''
(see \cite[Lemma 2.3]{GZ07} and compare with \cite[Lemma 3.5]{GZ07}).

\begin{lemma}  \label{lem:fdtl}
Let $\varphi, \psi \in {\mathcal E}^1(\Omega)$ be such that $\varphi \leq \psi$. Then
  $$
  E_\chi(\psi) \leq 2^{n+1} E_\chi(\varphi).
  $$
\end{lemma}

\begin{proof} 
Assume first $\varphi, \psi \in \mathcal T(\Omega)$ and observe that  
$
\int_\Omega h(-\varphi\slash  E_\chi(\varphi)) MA(\varphi)=1.
$
Fix  $\lambda > 0$ and recall that
$$
 \int_\Omega h(-\psi\slash 2 \lambda) d \mu_\psi = \int_0^{+\infty} h'(s) \int_{\{\psi \slash 2 \lambda < - s \}} (dd^c \psi)^n.
$$
Since $\varphi \leq \psi \leq 0$ we obtain
$$
\{\psi \slash 2 < - \lambda s\} \subset  \{\varphi \slash 2 < - \lambda s\} \subset \{ \varphi  < \psi \slash 2  - \lambda s \}
\subset \{ \varphi  <   - \lambda s \}.
$$ 
The comparison principle (Theorem \ref{thm:comparisonprinciple}) ensures that 
$$
\int_{\{\psi \slash 2 < - \lambda  s \}} MA(\psi\slash 2) \leq \int_{\{ \varphi < \psi \slash 2 - \lambda s\}} MA(\psi\slash 2) \leq  \int_{\{\varphi < - \lambda s\}} MA(\varphi).
$$
Hence 
$
 \int_\Omega h(-\psi\slash 2 \lambda) MA(\psi\slash 2)  \leq  \int_\Omega h(- \varphi \slash \lambda) MA(\varphi).
 $
We infer $E_\chi(\psi\slash 2) \leq E_\chi(\varphi)$ which implies  $E_\chi(\psi) \leq E_\chi(2 \varphi) \leq 2^{n+1} E_\chi(\varphi)$ by  the  Lemma \ref{lem:fdtl0}.

\smallskip

We now treat the general case. By Theorem \ref{thm:energy-approximation}, there exists a  sequence $({v}_j)$ in
 $\mathcal T(\Omega)$ decreasing to $\psi$ with $\lim_{j \to + \infty}  E_\chi(v_j) = E_\chi (\psi)$.
Pick any  sequence $(u_j)$ in $\mathcal T(\Omega)$ decreasing to $\varphi$.
 We set $f_j :=  \min \{u_j,v_j\}$ and consider
$$
 {\varphi}_j = P(u_j,v_j) := \sup \{ w ; w \in  PSH (\Omega) \text{ and }  w \leq f_j \}.
$$
Observe that $u_j +v_j  \in \mathcal T(\Omega)$ and $u_j +v_j \leq f_j \leq 0$,
hence   $u_j +v_j \leq  {\varphi}_j \leq v_j $. It follows that  ${\varphi}_j \in \mathcal T(\Omega)$. Since $(f_j)$ decreases to $\varphi$, it follows that $({\varphi}_j)$ decreases to $\varphi$.
The first step ensures that  for all $j \in \N$,
$$
E_\chi({\varphi}_j) \leq  2^{n+1} E_\chi(v_j)
$$
Using Theorem \ref{thm:energy-approximation} we infer
$$
 E_\chi({\varphi}) \leq \liminf_{j \to + \infty} E_\chi({\varphi}_j) \leq 2^{n+1} \liminf_{j \to + \infty} E_\chi(v_j) = 2^{n+1} E_\chi(\psi).
$$
\end{proof}

\begin{proposition}  \label{pro:highenergyequality}
A function $\varphi$ belongs to $\mathcal E_\chi(\Omega)$
if and only if $\varphi \in \mathcal E^1(\Omega)$ and $E_\chi(\varphi) < +\infty$.
Thus $\mathcal E_\chi(\Omega)$ is a positive cone.
Moreover $E_\chi(\varphi) = 0$ if and only if $\varphi = 0$.
\end{proposition}

We will show in Proposition \ref{pro:energycpct} that the cone $\mathcal E_\chi(\Omega)$ is furthermore convex.

\begin{proof}
Assume $\varphi \in \mathcal E_\chi(\Omega)$. 
By definition, there exists a sequence $(\varphi_j)$ in $\mathcal T(\Omega)$ decreasing to $\varphi$ such that $M:= \sup_j E_\chi(\varphi_j) < + \infty$. By lower semi-continuity Theorem \ref{thm:energy-approximation}, we have
$
E_\chi (\varphi) \leq \liminf_{j \to + \infty} E_\chi (\varphi_j) \leq M < + \infty.
$

\smallskip

Conversely if $\f \in \mathcal E^1(\Omega)$ and $E_\chi (\varphi)<+\infty$, we consider for $j \in \N$
$$
\f_j=\max( \f, j \rho) \in {\mathcal T}(\Omega).
$$
It follows from Lemma \ref{lem:fdtl} that $E_{\chi}(\f_j) \leq 2^{n+1} E_{\chi}(\f)$ is uniformly bounded,
hence $\varphi \in \mathcal E_\chi(\Omega)$, since $(\f_j)$ decreases to $\f$.

\smallskip

Fix $\gamma>0$. Recall that the Orlicz energy $E_\chi (\varphi)$ is finite if and only if there exists 
$\la>0$ such that $\int_{\Omega} h(-\f/\la) (dd^c \f)^n<+\infty$. 
 Since the Monge-Amp\`ere measure $(dd^c \f)^n$ is $n$-homogeneous, this finiteness condition holds
for $(\f,\la)$ if and only if it does so for $(\gamma \f, \gamma \la)$.
Thus $\mathcal E_{\chi}(\Omega)$ is a positive cone.

\smallskip

Clearly $\f=0$ implies $E_\chi(\varphi) = 0$. Assume now that 
$E_\chi(\varphi) = 0$.
 By definition there exists a sequence $(\lambda_j)$  of positive numbers decreasing to $0=E_\chi(\varphi)$ 
 such that for all $j \in \N$,
$$
\int_\Omega h(-\varphi\slash \lambda_j) MA(\varphi) \leq 1.
$$
Assume that $MA(\varphi) (\{\varphi < 0\}) > 0$. 
Then there exists a compact set $K \subset  \{\varphi < 0\}$ such that  $MA(\f)(K) > 0$. 
If we set $m_K = \max_K \varphi < 0$, then for all $j$ 
$$
h(-m_K \slash \lambda_j) \int_K MA(\varphi) \leq 1.
$$
Thus $E_\chi(\varphi) = \lim_j \lambda_j >  0,$ since $h(+ \infty) = + \infty$.
Therefore if $E_\chi(\varphi) = 0$ then $\varphi = 0$ almost everywhere for the measure $MA(\varphi)$,
and it follows from Lemma \ref{lem:classic} (extended to the class ${\mathcal E}^1(\Omega)$, see Remark \ref{rem:keylemE1})
that $\varphi \equiv 0$.
\end{proof}

\subsection{The subextension property}

Let $\widetilde \Omega \Subset \C^n$ be another bounded strongly pseudoconvex domain that contains 
$\overline{\Omega}$, hence $\Omega \Subset \widetilde \Omega \Subset \C^n$,
and fix $\varphi \in {\mathcal E}_\chi(\Omega)$.  
We consider  the maximal subextension of $\varphi$ to  $\widetilde \Omega$  defined by
 $$
 \widetilde \varphi := \sup \{u \in PSH^-(\widetilde \Omega) ; u \leq \varphi \, \, \text{in} \, \, \Omega\}.
 $$
It is not a priori clear that there exists a single plurisubharmonic function in 
$\widetilde \Omega$ that lies below $\f$; this follows from  some results established in \cite{GLZ19}:

 \begin{proposition} \label{pro:subext}
 Let $\varphi \in {\mathcal E}_\chi(\Omega)$. Then the following properties hold :
 \begin{enumerate}
 \item $\widetilde \varphi \in  {\mathcal E}_\chi(\widetilde \Omega)$ and $\widetilde \varphi  \leq  \varphi $ in $\Omega$;
 \item $(dd^c \widetilde \varphi)^n = 0$ in $\widetilde \Omega \setminus \overline \Omega$; 
 \item $(dd^c \widetilde \varphi)^n \leq  {\bf 1}_{\mathcal C} (dd^c \varphi)^n$ in $\Omega$, where 
 $\mathcal C :=\{ z \in \Omega ; \widetilde \varphi (z) = \varphi(z)\}$;
 \item $E_\chi(\widetilde \varphi) \leq E_\chi(\varphi).$
 \end{enumerate}
 \end{proposition}
 
 
 Similar properties have been established for the class $\mathcal F(\Omega)$  in \cite{CZ03,CH08}.
 
 \begin{proof}  
  Assume first that $\varphi \in \mathcal T(\Omega)$.  Since $\varphi$ is bounded from below in $\Omega$, it follows that
  $\widetilde \varphi \in PSH(\widetilde \Omega) \cap L^{\infty}(\widetilde \Omega)$ 
  and $\widetilde \varphi  \leq  \varphi $ in $\Omega$.
     We claim that $\widetilde \varphi \in \mathcal T(\widetilde \Omega)$. 
     Indeed, let $\widetilde \rho$ be a negative plurisubharmonic  exhaustion  for $\widetilde \Omega$. Fix $M >0$  such that $\varphi \geq - M$ in $\Omega$ and take a large constant $A > 1$ such that $A \widetilde \rho \leq - M $ in $\overline \Omega$. 
     Then $A \widetilde \rho \leq \widetilde \varphi $ in $\overline \Omega$. 
     This proves that $\widetilde \varphi \in \mathcal T (\widetilde \Omega)$.
     
     \smallskip
  
  A standard balayage method shows that $(dd^c  \widetilde \varphi)^n = 0$ in 
  $\widetilde \Omega \setminus \overline \Omega$,   proving $(2)$.
  Arguing as in the proof of \cite[Theorem 2.7]{GLZ19} one further shows that  $ (dd^c \widetilde \varphi)^n$ is carried by the contact set $\mathcal C :=  \{z \in \Omega ; \widetilde \varphi (z) =  \varphi(z)\}$.
  Indeed using  the notation of \cite{GLZ19}
   $\widetilde \varphi = P_{\widetilde \Omega} (h)$ is the psh envelope of $h = {\bf 1}_{\Omega} \varphi$ on $\widetilde \Omega$, which is quasi-lower semi-continuous in   $\widetilde \Omega$.
Since $\widetilde \varphi \leq  \varphi$ in $\Omega$,   Demailly's inequality yields
 \begin{equation}\label{eq:C}
 (dd^c \widetilde \varphi)^n  = {\bf 1}_{\mathcal C} (dd^c \widetilde \varphi)^n  \leq {\bf 1}_{\mathcal C} (dd^c  \varphi)^n,
 \end{equation}
 as Borel measures on $\Omega$, hence property $(3)$ in this case.
 
 \smallskip
 
Set $ \lambda := E_\chi( \varphi)$ so that $\int_\Omega h(- \varphi \slash \lambda) (dd^c  \varphi)^n = 1$.
Since $h(0) = 0$, It follows from the property $(2)$ and  \eqref{eq:C} that
$$
\int_{\widetilde \Omega} h(-\widetilde \varphi \slash \lambda) (dd^c \widetilde \varphi)^n 
\leq  \int_{\mathcal C} h(- \varphi \slash \lambda) (dd^c  \varphi)^n \leq 1.
$$
Hence $E_\chi( \varphi) = \lambda \geq E_\chi( \widetilde \varphi)$, proving $(4)$.

\medskip

We now treat the general case $\varphi \in E_\chi(\Omega)$. 
By definition, there is a decreasing sequence $(\varphi_j)$ in $\mathcal T (\Omega)$ 
converging to $\varphi$   and such that $B := \sup_j E_\chi(\varphi_j) < + \infty$.  
Let $\widetilde \varphi_j$ be the maximal subextension of $\varphi_j$ to $\widetilde \Omega$. 
The sequence $(\widetilde \varphi_j)$ decreases to $\widetilde \varphi$ in $\widetilde \Omega$,
and the first part of the proof yields $E_\chi( \widetilde \varphi_j) \leq E_\chi(\varphi_j) \leq B$ for any $j \in \N$. 
Thus $\widetilde \varphi \in \mathcal E_\chi(\widetilde \Omega)$ and  Proposition  \ref{pro:highenergyequality}
ensures that $E_\chi( \widetilde \varphi) \leq E_\chi(\varphi)$.

\smallskip

Taking limits in  $(dd^c \widetilde \varphi_j)^n  \leq {\bf 1}_{ \Omega} (dd^c  \varphi_j)^n$ 
 we obtain  $(dd^c \widetilde \varphi)^n  \leq {\bf 1}_{ \Omega} (dd^c  \varphi)^n$ in $\widetilde \Omega$.
 To check that the measure $(dd^c \widetilde \varphi)^n$ is carried by the contact set $\mathcal C$ 
 we  use the same argument as in the proof of \cite[Theorem 2.7]{GLZ19}.
 \end{proof}

 \section{Capacity characterization}
 
 In this section we provide an alternative definition of finite energy functions, estimating the 
 speed of decreasing of the capacity of their sublevel sets.

\subsection{Choquet integrals} 

Let $\mathrm{Cap}$ be the Monge-Ampère capacity
 with respect to $\Omega$ and  $h : \R^+ \longrightarrow \R^+$ a Young function.
  The integral with respect to $\mathrm{Cap}$ has been introduced in \cite{GSZ17}:
let $f : \Omega \longrightarrow \R^+$  be a positive Lebesgue measurable function;
the Choquet integral of $f$ with respect to  $\mathrm{Cap}$ and the weight $h$ is
 $$
 \int_\Omega h(f) d \mathrm{Cap} := 
 \int_0^{+\infty} h'(s) \mathrm{Cap}  (\{f > s\}) d s.
 $$

\begin{defi}
The Choquet-Orlicz space 
$L_h (\mathrm{Cap})$ is the set  of 
measurable functions 
$f : \Omega \longrightarrow \bar \R$ such that there exists $\lambda >0$ with $\int_\Omega h (\vert f\vert \slash \lambda) d \mathrm{Cap}< + \infty$. 

The Orlicz norm on $L_h ({\mathrm{Cap}})$  is defined by 
$$
 I_h (f) := \inf \left\{\lambda > 0 ; \int_\Omega h (\vert f\vert\slash \lambda) d \mathrm{Cap}  \leq 1 \right\}.
$$
\end{defi}


When $0<I_h (f)<+\infty$, the infimum is realized at $\la \in \R^+_*$ such that
$$
\int_\Omega h (\vert f\vert\slash \lambda) d \mathrm{Cap} = \int_0^{+\infty} \frac{1}{\la} h' \left( \frac{t}{\la} \right) \mathrm{Cap}(f>t) dt =1.
$$

\begin{proposition} \label{pro:choquet}
Let $f, g $ be two measurable functions in $\Omega$ and  $\alpha > 0$.
Then
\begin{enumerate}
\item $f \in L_h(\mathrm{Cap})  \Longleftrightarrow 0 \leq  I_h (f) < +\infty$;

\smallskip

\item $I_h(f) = 0$ if and only if $f = 0$;

\smallskip

\item if $f \in  L_h(\mathrm{Cap})$ and $0 \leq g \leq f $ then $g \in L_h(\mathrm{Cap}) $ and $I_h(g) \leq I_h(f)$;

\smallskip

\item if $(f_j)_{j \in \N}$ increases  to $f \in L_h(\mathrm{Cap})$,  then
$I_h(f) = \lim_{j \to + \infty} I_h(f_j)$;

\smallskip

\item 
if $f \in  L_h({\mathrm{Cap}})$, then $\alpha f \in L_h({\mathrm{Cap}})$ and 
 $I_h(\alpha f) = \alpha I_h(f)$;
\smallskip

\item if $f, g  \in  L_h(\mathrm{Cap})$ then 
$
 I_h (f + g) \leq   4 \max \left\{I_h (f), I_h (g)\right\}.
 $
\end{enumerate}
Thus $L_h(\mathrm{Cap})$ is a vector space.
\end{proposition}

\begin{proof}
1) By definition if $I_H (f) < \infty$ then $\int_\Omega h (\vert f\vert \slash \lambda) d \mathrm{Cap} \leq 1$
for any   $\lambda > I_h(f)$, hence $f \in  L_h (\mathrm{Cap})$.
Conversely if $f \in  L_h (\mathrm{Cap})$, there exists $\lambda_0 > 0$ such that $\int_\Omega h (\vert f\vert \slash \lambda_0) d \mathrm{Cap} < +\infty$. Since for any $s \geq 0$, $\lambda \mapsto h(s \slash \lambda)$ is decreasing in $\lambda $ and $\lim_{\lambda \to + \infty}h(s \slash \lambda) = 0$, it follows that 
$\lim_{\lambda \to + \infty} \int_\Omega h (\vert f\vert \slash \lambda) d \mathrm{Cap} = 0$. Hence there exists $\lambda > \lambda_0$ such that $\int_\Omega h (\vert f\vert \slash \lambda) d \mathrm{Cap} \leq 1$ so that 
$I_h(f) \leq \lambda < + \infty$.

\smallskip

2) If $f = 0$ a.e. in $\Omega$, then $h (\vert f\vert \slash \lambda) = 0$ 
a.e. in $\Omega$ since $h(0) = 0$,  hence $I_h(0) = 0$.
Assume conversely that $f \neq 0$.  By definition
$ \mathrm{Cap}(\{|f|>0\})  > 0$.
By the upper semi-continuity property of the capacity there exists $\tau > 0$ such that $ \mathrm{Cap}(F_\tau)  > 0$,
where $F_\tau := \{\vert f\vert > \tau\}$. 
Now
  \begin{eqnarray*} 
  \int_\Omega h (\vert f\vert\slash \lambda) d  \mathrm{Cap} &\geq &
 \int_0^{\tau} (1\slash \lambda) h'(s\slash \lambda)  \mathrm{Cap} (\{\vert f\vert \geq s\}) ds \\
 &\geq & \mathrm{Cap}(F_\tau)  \int_0^{\tau} (1\slash \lambda) h'(s\slash \lambda) d s 
 = h  (\tau\slash \lambda) \mathrm{Cap}(F_\tau).
\end{eqnarray*}
In particular $\mathrm{Cap}(F_\tau) \leq h (\tau\slash\lambda)  \int_\Omega h  (\vert f\vert \slash \lambda) d  \mathrm{Cap}$.
 Since $\lim_{\lambda \to 0} h ( (\tau\slash \lambda) = +\infty$
 we infer $ \lim_{\lambda \to 0} \int_\Omega h(\vert f\vert \slash \lambda) d  \mathrm{Cap} = +\infty$, hence $I_h(f) > 0$.
 
 \smallskip
 
3) is clear, we now treat 4).  Set $\la_j=I_h (f_j)$, so that
$$
\int_0^{+\infty} \frac{1}{\la_j} h' \left( \frac{t}{\la_j} \right) \mathrm{Cap}(f_j>t) dt =1.
$$
The sequence $(\la_j)$ is increasing since so is $j \mapsto f_j$. Set $\la=\lim_j \la_j$. 
The continuity properties of the Monge-Amp\`ere capacity $\mathrm{Cap}$ ensure that 
$$
\int_0^{+\infty} \frac{1}{\la} h' \left( \frac{t}{\la_j} \right) \mathrm{Cap}(f>t) dt =1,
$$
hence $\la=I_h (f)=\lim I_h (f_j)$ as desired.

\smallskip
 
5) Fix $\lambda > 0$ and $\alpha > 0$ and assume $f \geq 0$. 
Setting $ s := \alpha t$ we obtain
\begin{eqnarray*}
\int_\Omega h (\alpha \vert f\vert \slash \lambda) d  \mathrm{Cap} 
&=& \int_0^{+ \infty} (1 \slash \lambda) h'(s\slash \lambda) \mathrm{Cap} (\{\alpha f \geq s\}) d s \\
&= & \int_0^{+ \infty} (1 \slash \lambda) h'(\alpha t\slash \lambda)  \mathrm{Cap} (\{ f \geq t \})  \alpha d t \\
&=&  \int_0^{+ \infty} (\alpha \slash \lambda) h'(\alpha t\slash \lambda)  \mathrm{Cap} (\{ f \geq t \})  d t.
\end{eqnarray*}
If $\lambda = \alpha I_h (f)$, then the last integral in the previous inequalities is equal to $1$.
Hence $\lambda = I_h(\alpha f)$ and   $I_h( \alpha f) = \alpha I_h(f)$.

\smallskip

We finally prove  6).  The subadditivity of the capacity yields the following quasi-triangle inequality : 
for any positive measurable functions $u, v$,
$$
\int_\Omega (u + v) \, d \, \mathrm{Cap} \leq 2 \left(\int_\Omega u \, d \, \mathrm{Cap} + \int_\Omega v \, d \, \mathrm{Cap}\right).
$$
Assume $f, g \geq 0$ and fix $\lambda = \max\{I_h(f), I_h(g)\}$. 
Using the convexity of $h$ and the quasi-triangle inequality, we obtain
\begin{eqnarray*}
\int_\Omega h ((f + g)\slash 4 \lambda)  d  \mathrm{Cap} &\leq &  \int_\Omega (1\slash 2) h (f \slash \lambda)  d  \mathrm{Cap} +  \int_\Omega (1\slash 2) h (g \slash \lambda)  d  \mathrm{Cap} \leq 1.
\end{eqnarray*}
Hence $4 \lambda \geq I_h(f+g)$ which proves the required inequality.
\end{proof}

The following convexity property will play an important role in Section \ref{sec:ma}.

\begin{proposition} \label{pro:choquetinfini}
  Let $(f_j)$ be a sequence of positive measurable functions in $\Omega$ such that  $(I_h(f_j))$ is bounded.
Fix $(\alpha_j)_{j \in \N}, (\varepsilon_j)_{j \in \N}$  sequences of positive numbers such that 
$\sum_{j \in \N}  \alpha_j, \sum_{j \in \N}  \e_j \leq 1$. 
  Then 
$
 I_h \left(\sum_{j \in \N}  \varepsilon_j \alpha_j f_j\right)  \leq \sup_j I_h(f_j).
$
\end{proposition}

\begin{proof}
Observe that  if $(u_j)_{j\in \N}$ is a sequence of measurable functions in $\Omega$, then
$$
\left \{\sum_{j \in \N} \varepsilon_j u_j \geq  t \right \} \subset \bigcup_{j\in \N} \left\{ u_j \geq t\right\}
$$
The subadditivity of the capacity yields 
\begin{equation} \label{eq:Subineq}
\int_\Omega\left(\sum_{j \in \N} \varepsilon_j u_j\right) d  \mathrm{Cap} \leq \sum_{j \in \N}  \int_\Omega u_j d  \mathrm{Cap}.
\end{equation}
 Fix $\lambda > 0$ and set $f = \sum_{j \in \N} \varepsilon_j \alpha_j f_j$. 
We obtain
 $
 h (f\slash \lambda) \leq \sum_{j \in \N}  \varepsilon_j h(\alpha_j f_j\slash \lambda)
 $
 by convexity of $h$.
Together with \eqref{eq:Subineq},  this yields
  $$
 \int_\Omega h (f\slash \lambda) d  \mathrm{Cap}  
 \leq \sum_{j \in \N}  \int_\Omega  h \left(\alpha_j f_j \slash \lambda\right) d  \mathrm{Cap} 
 \leq \sum_{j \in \N}  \alpha_j \int_\Omega  h \left( f_j \slash \lambda\right) d  \mathrm{Cap}.
 $$
 
Set $\lambda = \sup_{j \in \N}  I_h(f_j)$. Then
 $
 \int_\Omega  h \left( f_j \slash \lambda\right) d  \mathrm{Cap} 
 \leq  \int_\Omega  h \left( f_j \slash I_h(f_j)\right) d  \mathrm{Cap} 
 = 1,
 $
 hence $\int_\Omega h (f\slash \lambda) d  \mathrm{Cap}  \leq  \sum_{j \in \N} \varepsilon_j \leq 1$ 
 which yields $\lambda  \geq  I_h  (f)$.
 \end{proof}

\subsection{Alternative definition of high energy classes}

Fix $\chi \in {\mathcal W}^+$ and $h$ the associated Young function.
We consider 
 $
 \widetilde h(s) := \int_0^s t^n h'(t)  dt
 $
 and set 
$$
J_\chi(\varphi) := I_{\widetilde h} (-\varphi).
$$

\begin{defi}
The class $\widetilde{\mathcal E}_\chi (\Omega)$ is the set of functions $u \in PSH(\Omega)$ 
such that there exists   $(u_j)$ in $\mathcal T (\Omega)^{\N}$ decreasing to $u$ with
$\sup_j J_{\chi}(-u_j) < + \infty$.
\end{defi}

 By definition $\widetilde{\mathcal E}_\chi (\Omega) \subset L_{\widetilde h} ({\mathrm{Cap}})$ is a  convex cone.
We   show that $\mathcal E_\chi (\Omega)=\widetilde{\mathcal E}_\chi (\Omega)$.

\begin{proposition} \label{prop:cap-characterization}
\text{ }

\begin{enumerate}
\item    $\varphi \in \mathcal E_\chi (\Omega)$ if and only if  $\varphi \in \widetilde{\mathcal E}_\chi (\Omega)$.
Moreover
\begin{equation*}
J_\chi(\varphi) \leq 2 \max \{1, E_\chi (\varphi)\} \, \, \,  \text{and} \, \, \, E_\chi(\varphi) \leq \max \{1, J_\chi (\varphi)^{n+1}\}.
\end{equation*}

\item  If $(u_j)_{j \in \N}$ are plurisubharmonic functions  decreasing  to $u \in \widetilde{\mathcal E}_\chi (\Omega)$ ,  then 
$$
J_\chi(u) = \lim_{j \to + \infty} J_\chi(u_j).
$$ 
 In particular the sequence $(E_\chi(u_j))_{j \in \N}$ is bounded.
\item If $\varphi \in \mathcal E_\chi (\Omega)$ and $\la=J_\chi(\varphi)$,  then
 for all $s > 0$
 $$
 \mathrm{Cap} (\{\varphi \leq - s \}) \leq \frac{1}{\widetilde h(s \slash \lambda)}\cdot
 $$ 
Conversely if there exists $\la_0>0$ and $C_0>0$ such that for all $s > 0$
 $$
 \mathrm{Cap} (\{\varphi \leq - s \}) \leq \frac{C_{0}}{(1+s)^{2} s^n h'(s/\lambda_0)},
 $$ 
then $\f \in \mathcal E_\chi (\Omega)$ and $E_\chi(\varphi) \leq \max \{C_0, \lambda_0\}^{n+1}$.
\end{enumerate} 
\end{proposition}

The last information underlines the interest of the capacity point of view: 
we now know that a function belongs to a high energy class
if and only if the capacity of its sublevel sets decreases rapidly enough,
while it is usually tricky to extract information from the Monge-Amp\`ere measure.

\begin{proof} 
Assume first that $\varphi \in \mathcal T (\Omega)$ and fix $\lambda \geq E_\chi (\varphi)$. Observe that
$$
 \int_\Omega \widetilde h (-\varphi\slash \lambda) d \mu_\varphi = \int_0^{+\infty} (1\slash \lambda^{n+1}) h'(s\slash \lambda) 
   s^n MA(\f)(\{\varphi  < - s \}) ds \leq 1.
$$
Applying Lemma \ref{lem:classic} with $s > 0$ and $t =  s$   we obtain 
\begin{eqnarray*}
\int_\Omega \widetilde h (-\varphi \slash 2 \lambda) d \mathrm{Cap}  &=&\int_0^{+\infty}  (1\slash \lambda^{n+1}) h'(s\slash \lambda)  s^n   \mathrm{Cap} (\{\varphi < - 2 s  \}) d s \\
&\leq& (1 \slash \lambda^{n}) \int_\Omega h(-\varphi \slash \lambda)  MA(\varphi) \\
&\leq& \int_0^{+\infty}  (1\slash \lambda^{n+1}) h'(s\slash \lambda)  s^n   \mathrm{Cap} (\{\varphi < - s  \}) d s \\
&=& \int_\Omega \widetilde h (-\varphi \slash \lambda) d \mathrm{Cap}.
\end{eqnarray*}
This implies that for any $\lambda > E_\chi (\varphi)$, 

$$
\int_\Omega \widetilde h (-\varphi \slash 2 \lambda) d \mathrm{Cap}   \leq (1 \slash \lambda^{n}) \int_\Omega  h(-\varphi \slash \lambda)  d \mu_\varphi \leq  \int_\Omega \widetilde h (-\varphi \slash \lambda) d \mathrm{Cap}.
$$
Hence if $\lambda > \max \{1, E_\chi (\varphi)\}$ then $\int_\Omega \widetilde h (-\varphi \slash 2 \lambda) d \mathrm{Cap}   \leq  1$, i.e. $\lambda \geq J_\chi(\varphi \slash 2)$. Which implies that $J_\chi(\varphi) \leq 2 \lambda \leq 2 \max \{1, E_\chi (\varphi)\}$.

When $\varphi$ merely belongs to $\mathcal E_\chi(\Omega)$, 
we use the approximating sequence $(\widetilde{\varphi}_j)$ constructed in Theorem \ref{thm:energy-approximation}
 and apply the previous inequality. 
 The conclusion follows by taking the limit and using  Proposition \ref{pro:choquet}.

\medskip

On the other hand if $\lambda >  \max \{1, J_\chi (\varphi)\}$, it follows from the convexity of $h$ that
$$
\int_\Omega  h(-\varphi \slash \lambda^{n+1})  d \mu_\varphi \leq (1 \slash \lambda^{n}) \int_\Omega  h(-\varphi \slash \lambda)  d \mu_\varphi \leq  \int_\Omega \widetilde h (-\varphi \slash \lambda) d \mathrm{Cap} \leq 1.
$$
Hence $\lambda^{n+1} \geq E_\chi(\varphi)$, which yields $E_\chi(\varphi) \leq \max \{1, J_\chi (\varphi)^{n+1}\}$.
Thus
\begin{equation} \label{eq:Energy-capacity} 
 J_\chi(\varphi \slash 2) \leq \max\{1, E_\chi(\varphi) \}  \, \, \, \text{and} \, \, \, E_\chi(\varphi) \leq \max \{1, J_\chi(\varphi)^{n+1}\}.
\end{equation}

\smallskip

Conversely if   $(\varphi_j)$ is any  sequence  in $\mathcal T(\Omega)$ 
decreasing to $\varphi \in \widetilde {\mathcal E}_\chi(\Omega)$,
 it follows from \eqref{eq:Energy-capacity} that the sequence $E_\chi(\varphi_j)$ is bounded.
Fix $\lambda > \liminf_{j \to + \infty} E_\chi (\varphi_j)$. 
Along an appropriate subsequence  we have 
$\lambda >  E_\chi (\varphi_j)$ and 
$$
\int_\Omega h(-\varphi_j \slash  \lambda)MA(\varphi_j)  \leq \int_\Omega h(-\varphi_j \slash  E_\chi (\varphi_j)) MA(\varphi_j)\leq 1.
$$
By lower semi-continuity we obtain
$$
\int_\Omega h(-\varphi \slash  \lambda) MA(\varphi) \leq \liminf_{j \to + \infty} \int_\Omega h(-\varphi_j \slash  \lambda)MA(\varphi_j)  \leq 1.
$$
Hence $\lambda \geq E_\chi(\varphi)$, showing that $E_\chi (\varphi) \leq \liminf_{j \to + \infty} E_\chi (\varphi_j).$
We infer
$$
E_\chi (\varphi) \leq \liminf_{j \to + \infty} E_\chi (\varphi_j) \leq  \max \{1, J_\chi  (\varphi)^{n+1}\}.
$$
Thus $\varphi \in \mathcal E_\chi(\Omega)$ as claimed.

\medskip

The continuity property (2) follows from Proposition \ref{pro:choquet}.4. We now prove (3).
It follows from Chebyshev inequality that if $\f \in L_{\widetilde h}({\rm Cap})$ then for all $s > 0$
 $
 \mathrm{Cap} (\{\varphi \leq - s \}) \leq \frac{1}{\tilde h(s \slash \lambda)},
 $
 where $\la=J_{\chi}(\f)$.

Conversely if there exists $\la_0, C_0 > 0$ such that 
 $
 \mathrm{Cap} (\{\varphi \leq - s \}) \leq \frac{C_{0}}{(1+s)^{2} s^n h'(s/ \lambda_0)},
 $
then for $\lambda > 0$ we obtain
$$
\int_{\Omega} \widetilde h(|\f| \slash \lambda) d{\rm Cap} \leq  \int_0^{+ \infty}  \frac{C_0 s^n h'(s\slash \lambda)}{\lambda (1+ s)^{2} s^n  h'(s\slash \lambda_0)} ds.
$$
Thus if $\lambda \geq \max \{\lambda_0, C_0\} $, the last integral is bounded from above by $1$
hence $J_{\chi}(\f) \leq \max \{\lambda_0, C_0\}$,
and  $\f \in \mathcal E_\chi (\Omega)$ with $E_\chi (\varphi) \leq \max\{\lambda_0,C_0\}^{n+1}$.
\end{proof}

\begin{example} 
Assume that $h(t) := t^p \slash p$ with $p \geq 1$. Then $\tilde h(t) = t^{n + p} \slash (n+p)$ and for any $\varphi \in \mathcal E_\chi(\Omega)$, we obtain
$$
E_\chi (\varphi) = \left(\int_\Omega (-\varphi)^p MA(\varphi)\right)^{1 \slash p}, \, \, J_\chi (\varphi) = \left(\int_\Omega (-\varphi)^{n + p} d{\rm Cap}\right)^{1\slash n + p}.
$$
The identity $\widetilde {\mathcal E}_p (\Omega) = \mathcal E_p(\Omega)$ has been established in \cite[Proposition B]{BGZ09}.
\end{example}

\begin{exa}
When $\Omega=\B$ is the unit ball of $\C^n$ and  
the functions under consideration are {\it radial}, one can easily determine whether they
belong to a class ${\mathcal E}_{\chi}(\B)$.
Consider e.g. $L_k =\log \circ \cdots \circ \log$ the $k^{th}$-iterated logarithm, and
$$
\f_k(z)=-L_k(-\log ||z||+C_k),
$$
where $C_1=1$ and $C_{k+1}=\exp(C_k)$ so that ${\f_k}_{|\partial \B}=0$.

Recall \cite[Ex. 4.37]{GZbook} that ${\rm Cap}(\B(e^{-t})=t^{-n}$, so
$
{\rm Cap}(\f_k <-t) \leq \frac{C_k'}{e_k(t)^n},
$
where $e_k=\exp \circ \cdots \circ \exp$ is the $k^{th}$-iterated exponential.
It therefore follows from Proposition \ref{prop:cap-characterization} that 
$\f_k \in {\mathcal E}_{\chi_k}(\B)$, where $\chi_k(t)=-e_k(-t)$.
\end{exa}
 
 One can also easily produce {\it toric examples}, assuming
 $\Omega$ merely has $(S^1)^n$-symmetries. We refer the interested reader
 to  \cite{CGSZ19} for a discussion of toric finite energy classes
 in the context of compact K\"ahler manifolds.

 \subsection{Moser-Trudinger inequalities}

As   recalled in Section \ref{sec:prelim2},  plurisubharmonic functions have strong
integrability properties. For functions $\f$ in a finite energy class
${\mathcal E}_{\chi}(\Omega)$, the fastest the growth of $\chi$ at infinity, the stronger
the integrability properties of $\f$.
Let $h: \R^+ \rightarrow \R^+$ be a Young function, and
let $\chi(t)=-h(-t) \in {\mathcal W}^+$ denote the weight associated to $h$. We consider
$$
\tilde{h}: t \in \R^+ \mapsto \int_0^t s^n h'(s) ds \in \R^+,
$$
and the Young function $H_{\beta}$  defined by
$
H_{\beta}(t)= \exp \left( \beta [\tilde{h}(t)]^{1/n} \right).
$

\begin{theorem} \label{thm:mosert}
Fix $\beta<2n$.
 There is $C(\beta,{\chi})>0$ such that for all
$\f \in {\mathcal E}_{\chi}(\Omega)$, 
$$
\int_{\Omega} H_{\beta} \left(\frac{-\f}{2 \max(1,E_{\chi}(\f))} \right) dV \leq C(\beta, {\chi}).
$$
\end{theorem}

This Moser-Trudinger inequality is a generalization of \cite[Theorem 4.6]{DGL21},
which dealt with polynomial weights. 
We refer the reader to the comments after \cite[Theorem 4.6]{DGL21} for an
historical account of this important inequality.
If the growth of $h$ is at least exponential,
then $\log \tilde{h}(t) \sim \log h(t)$ hence -absorbing error terms  by varying $\beta$-
the result holds with $H_{\beta}(t)=\exp \left( \beta [{h}(t)]^{1/n} \right)$.

 
\begin{proof}
   It follows from \cite[Proposition 6.1]{ACKPZ09} 
that for all $\beta<2n$,  there exists $C_{\beta}>0$ such that
  $$
 {\rm Vol}(\{\varphi <-t\}) \leq C_{\beta} \exp\left( -\frac{\beta}{{\rm Cap}(\{\varphi <-t\})^{1/n}} \right)
 \leq C_{\beta} \exp\left( - \beta \tilde{h}(t \slash \lambda)^{1/n} \right),
  $$
  where the last inequality follows  from Proposition \ref{prop:cap-characterization} with $\lambda = J_\chi(\varphi)$.
We infer
  $$
  {\rm Vol}(\{\varphi \slash \lambda <-t\}) \leq C_{\beta} \exp\left( - \beta \tilde{h}(t)^{1/n} \right),
  $$
  
  Let $\sigma$ be the Young function defined by 
  $
  \sigma'(t)= \frac{1}{1+t^2} \exp \left( \beta [\tilde{h}(t)]^{1/n} \right).
  $
 It follows that
  \begin{eqnarray*}
  \int_{\Omega} \sigma (-\f\slash \lambda) dV & =& 
  \int_0^{+\infty} \sigma'(t) {\rm Vol} (\f<-t) dt \\
  & \leq &  C_{\beta} \int_0^{+\infty}   \frac{dt}{1+t^2}=C_{\beta} <+\infty
  \end{eqnarray*}
  
Since $h(t) \geq t$, the term $\frac{1}{1+t^2}$ in the definition of 
$\sigma'(t)$ is negligible by comparison to the exponential growth of 
$\exp \left( \beta [\tilde{h}(t)]^{1/n} \right)$. It can be absorbed by adjusting the value of $\beta$.
Similarly $\log \sigma'(t)$ and $\log \sigma(t)$ have equivalent growth at infinity.
Thus varying $\beta$ -hence adjusting $C_{\beta}$- we conclude that
$$
\int_{\Omega} H_{\beta}(-\f\slash J_\chi(\varphi)) \, dV \leq C'(\beta, {\chi}).
$$
Since $J_\chi(\varphi) \leq 2 \max \{1,E_\chi(\varphi)\}$ the required inequality follows.
\end{proof}

\section{Range of the Monge-Ampère operator} \label{sec:ma}

 In this section we study the range of the complex Monge-Ampère operator,
extending Theorem \ref{thm:ceg} to weights that have arbitrary growth at infinity.

 \subsection{Compactness properties}
 
 Functions with uniformly bounded energies are relatively compact for the induced  $L^1$-topology.

 \begin{prop} \label{pro:energycpct}
 Fix $A>0$. The set
 $$
 {\mathcal E}_{\chi}^A(\Omega):=\left\{ \f \in {\mathcal E}_{\chi}(\Omega), \; E_{\chi}(\f) \leq A \right\}
 $$
 is a relatively compact subset of ${\mathcal E}_{\chi}(\Omega)$. 
 Moreover  if $(\f_j) \in {\mathcal E}_{\chi}^A(\Omega)$ then
 $$
 \p=\sum_{j \geq 1} 4^{-j} \f_j \in {\mathcal E}_{\chi}^{B} (\Omega),
 $$
 with $B=(2A)^{n+1}$.
In particular ${\mathcal E}_{\chi}(\Omega)$ is a convex cone.
 \end{prop}

  \begin{proof}
 Fix $\f_j \in {\mathcal E}_{\chi}^A(\Omega)$ and set $\la_j=E_{\chi}(\f_j) \leq A$.
 Extracting and relabelling we can assume that $\la_j \rightarrow \la \leq A$. Since $0 \leq t \leq h(t)$ we obtain
 $$
 \int_{\Omega} (-\f_j) MA(\f_j) \leq \la_j \int_{\Omega} \left(\frac{-\f_j}{\la_j} \right) MA(\f_j)
 \leq \la_j \int_{\Omega} h \left(\frac{-\f_j}{\la_j} \right) MA(\f_j) =\la_j \leq A.
 $$
 
 The sequence $(\f_j)$ is therefore relatively compact in the Cegrell class ${\mathcal E}^1(\Omega)$.
 Extracting and relabelling we can   assume that $\f_j \rightarrow \f \in {\mathcal E}^1(\Omega)$ a.e. and in $L^1$, hence
 $
 \tilde{\f}_j=\left( \sup_{\ell \geq j} \f_{\ell} \right)^*
 $
 decreases pointwise to $\f$.   Lemma \ref{lem:fdtl} ensures that
 $$
 E_{\chi}(\tilde{\f}_j) \leq 2^{n+1} E_{\chi}(\f_j) \leq 2^{n+1} A,
 $$
 hence $\f \in  {\mathcal E}_{\chi}^{A'}(\Omega) \subset  {\mathcal E}_{\chi}(\Omega)$,
 with $A'=2^{n+1} A$.
 
 \smallskip
 
 Fix now $\f_j \in {\mathcal E}_{\chi}^A(\Omega)$  and consider
 $ \p=\sum_{j \geq 1} 4^{-j} \f_j$.
 We can assume without loss of generality that
 $J_{\chi}(\f_j),E_{\chi}(\f_j) \geq 1$.
 It follows from Proposition \ref{pro:choquetinfini} 
 (with $\e_j=\alpha_j=2^{-j}$) and Proposition \ref{prop:cap-characterization} that
 $$
 E_{\chi}(\p) \leq \left( 2A \sum_{j \geq 1} \e_j \sum_{j \geq 1} \alpha_j \right)^{n+1}
 =(2A)^{n+1}.
 $$
     This shows in particular that ${\mathcal E}_{\chi}(\Omega)$ is a convex cone.
  \end{proof}

 \subsection{Integrability}
 
 In this section we analyze the  condition $\mathcal E_\chi(\Omega) \subset L_{\chi}(\mu)$, both from
 a qualitative and a quantitative point of view.
 While the qualitative integrability condition is good enough for polynomial weights, we will see in Example \ref{exa:n-Entropie} that
 it is too weak for weights with exponential growth.

 \subsubsection{Qualitative integrability}

 \begin{defi}
 Fix $\chi \in {\mathcal W}^+$ and $a>0$. 
 A positive Radon measure $\mu$ satisfies property $I(a,\chi)$ if there exists $C_a>0$ such that 
 for all $\f \in \mathcal E_\chi(\Omega)$,
 \begin{equation}  \label{eq:Ia}
 \tag{$I(a,\chi)$}
 ||\f||_{L_{\chi}(\mu)} \leq a E_{\chi}(\f)+C_a.
 \end{equation}
 \end{defi}
 
 This inequality is equivalent to the integrability property $\mathcal E_\chi(\Omega) \subset L_{\chi}(\mu)$:
 
 \begin{prop} \label{pro:qualitative}
 Fix $\chi \in {\mathcal W}^+$.
 The following properties are equivalent:
 \begin{enumerate}
 \item $\mathcal E_\chi(\Omega) \subset L_{\chi}(\mu)$;
 \item there exists $a>0$ such that the measure $\mu$ satisfies $I(a,\chi)$.
 \end{enumerate}
 \end{prop}

  \begin{proof}
 The implication $(2) \Rightarrow (1)$ is clear so it suffices to show that 
 if  $\mu$ does not satisfy any of the properties $I(a,\chi)$, then there exists 
 $\p \in \mathcal E_\chi(\Omega)$ such that $||\p||_{L_{\chi}(\mu)}=+\infty$.
Observe that  $\mu$ satisfies $I(a,\chi)$ for some $a>0$ if and only if
  $$
  ||\f||_{L_{\chi}(\mu)} \leq \max(a,C_a) \left[ E_{\chi}(\f)+1 \right].
  $$
  Assume this is not the case.
 Extracting and relabelling if necessary, we can  find 
 $\f_j \in \mathcal E_\chi(\Omega) $ such that 
 $$
 ||\f_j||_{L_{\chi}(\mu)} \geq 8^j \left[ E_{\chi}(\f_j)+1 \right].
 $$
 
 If $E_{\chi}(\f_j) \leq M$ is uniformly bounded, it follows from 
 Proposition \ref{pro:energycpct} that 
 $\p=\sum_j 4^{-j} \f_j \in \mathcal E_\chi(\Omega)$. 
 Since $\p \leq  4^{-j} \f_j \leq 0$ for all $j \in \N$, we infer
 $$
 ||\p||_{L_{\chi}(\mu)} \geq ||4^{-j} \f_j||_{L_{\chi}(\mu)}=4^{-j} || \f_j||_{L_{\chi}(\mu)}  \geq 2^j,
 $$
 hence $||\p||_{L_{\chi}(\mu)}=+\infty$.
 
 The same reasoning can be applied if $(E_{\chi}(\f_j))$ admits a bounded subsequence.
 We finally assume that $1 \leq M_j=E_{\chi}(\f_j) \rightarrow +\infty$. Consider
 $\p=\sum_j 4^{-j} \frac{\f_j}{M_j}$.
 It follows from Proposition \ref{pro:choquet}.4 and 
Proposition \ref{pro:energycpct}
 that 
 $$
 E_{\chi}\left( \frac{\f_j}{M_j} \right) \leq \frac{1}{M_j} (E_{\chi}(\f_j)) \leq 1,
 $$
 hence $\p  \in \mathcal E_\chi(\Omega)$. Now
 $
 ||\p||_{L_{\chi}(\mu)} \geq \frac{4^{-j}}{M_j} || \f_j||_{L_{\chi}(\mu)}  \geq 2^j,
 $
so $||\p||_{L_{\chi}(\mu)}=+\infty$.
 \end{proof}

 We exhibit a large family of measures satisfying this  integrability property.
 
 \begin{exa} \label{exa:n-Entropie}
 Assume $\mu=f dV$ 
 is such that  $ \int_\Omega f (\log (1+ f))^p dV < + \infty$, for some $p>0$.
 We claim that  $\mathcal E_\chi (\Omega) \subset L_\chi(\Omega, \mu)$ for any  weight $\chi \in \mathcal W^+$
 which has at least exponential growth.
 
 Indeed let $h(s) := -\chi(-s)$ be the Young function associated to $\chi \in {\mathcal W}^+$, 
 and fix $\varphi \in \mathcal E_\chi (\Omega)$. We assume without loss of generality that
 $E_{\chi}(\f)=1$.  It follows from Theorem \ref{thm:mosert} that  
 $\int_\Omega \exp ({\beta  h(-\varphi)^{1\slash n}}) d V < + \infty$, for all $0 < \beta < 2 n $.

 We need to show that   $\int_\Omega h(-\varepsilon \varphi) f d V < +\infty$ for some $\varepsilon >0$. Since
 $$
 \int_\Omega h(-\varepsilon \varphi) f d V \leq \int_\Omega \theta (f) d V + \int_\Omega \theta^* \circ h(-\varepsilon \varphi) dV,
 $$
 where $\theta^*$ is the Legendre transform of   function $t \longmapsto \theta (t) = t [\log (1+t)]^p$,
 it suffices to check that $ \int_\Omega \theta^* \circ h(-\varepsilon \varphi) dV < +\infty$ for some $\e>0$.
  We let the reader check that $\theta^*(s) \sim s^{1-1/ p} \exp({s^{1/ p}})$  as  $s \to + \infty$.
 Adjusting the value of $\beta$, it is thus enough to show that there exists $0 <\varepsilon < 1$ small enough so that for $s \geq 1$,
 $$
h(\varepsilon s)  \leq    \beta h(s)^{p/n} +C \leq C_2 h(s)^{p/n}.
 $$

Since $h$ has at least exponential growth, we can assume that $h \sim \exp(H) $ with $H$ convex.
The previous inequality can thus be written as
$
H(\e s) \leq \frac{p}{n} H(s)+C_3.
$
Choosing $\e=p/n$, the latter holds true by convexity, up to enlarging $C_3$.
\end{exa}

  \begin{remark}
If $\mu=f dV$  has finite $p$-entropy $ \int_\Omega f (\log (1+ f))^p dV < + \infty$, 
then $\mu \leq C \left( {\rm Cap} \right)^{p/n}$ is dominated by the Monge-Amp\`ere capacity
(see \cite[Lemma 4.2]{Kol05}).
It follows that $\mu=(dd^c \f)^n$ for some unique $\f \in {\mathcal E}^1(\Omega)$, which moreover
\begin{itemize}
\item belongs to ${\mathcal E}^q(\Omega)$ for $q<\frac{pn}{n-p}$ if $0<p<n$ (see \cite[Proposition 5.3]{BGZ09});
\item belongs to ${\mathcal E}_{exp}(\Omega)$  if $p=n$ (see \cite[Theorem 5.1]{BGZ09});
\item is continuous if $p>n$ (see \cite[Theorem 4.6]{Kol05}).
\end{itemize}
Examples of radial functions show that these bounds are essentially optimal.
 \end{remark}

 \subsubsection{Strong coercivity}
 
 We now show that if $\mu$ is the Monge-Amp\`ere measure of a function in $\mathcal E_\chi(\Omega)$,
 then $\mathcal E_\chi(\Omega) \subset L_{\chi}(\mu)$ and we can moreover reinforce
 the  coercivity inequality provided by Proposition \ref{pro:qualitative}.

\begin{prop} \label{pro:quantitative}
Fix $\chi \in {\mathcal W}^+$, $\psi \in \mathcal E_\chi(\Omega)$ and set $\mu = MA(\psi)$.
The measure $\mu$ satisfies $I(a,\chi)$ for all $a>1$.
\end{prop}

\begin{proof}
It follows from Theorem \ref{thm:energy-approximation} that 
it suffices to check $I(a,\chi)$ 
for tests $\f$ that belong to  $\mathcal T(\Omega)$.
Assume first also
that $\psi \in \mathcal T(\Omega)$. Set  $ \varepsilon := a -1 > 0$ and let $c > 0$ and $\lambda > 0$ be fixed. Then 
$$
\int_\Omega h(-\varphi \slash a \lambda) MA(\psi) = \int_0^{+ \infty} h'(s) \int_{\{\varphi \slash \lambda < - a s\} } MA(\psi).
$$
Observe that 
$$
\{\varphi \slash \lambda < - (1+ \varepsilon) s\} \subset \{\varphi   < c \psi  -  \lambda s\} \cup \{ c \psi   <  -  \lambda \varepsilon s\}\cdot
$$
The comparison principle (Theorem \ref{thm:comparisonprinciple}) yields
$$
\int_{\{\varphi   <c \psi   -  \lambda s\} } MA(\psi) \leq c^{-n} \int_{\{\varphi <c \psi  - \lambda s\} }  MA(\varphi) \leq c^{-n} \int_{\{\varphi  \slash \lambda < -  s\} }  MA(\varphi),
$$
hence 
$$
\int_\Omega h (-\varphi \slash a \lambda) d \mu_\psi \leq  c^{-n}  \int_\Omega h(-\varphi \slash \lambda)d  \mu_\varphi + 
 \int_{\Omega } h(-c \psi\slash \varepsilon \lambda) d \mu_\psi.
$$

We choose  $c=2^{1/n}$ and fix $\lambda \geq E_\chi(\varphi)$ so that
$$
\int_\Omega h (-\varphi \slash a \lambda) d \mu_\psi \leq  \frac{1}{2} + 
 \int_{\Omega } h(-c \psi\slash \varepsilon \lambda) d \mu_\psi.
$$
If we moreover impose  $\lambda \varepsilon \slash c \geq 2 E_\chi(\psi)$, the convexity of $h$ yields
$$
\int_{\Omega } h(-c \psi\slash \varepsilon \lambda) d \mu_\psi \leq \int_{\Omega } h(- \psi\slash 2  E_\chi(\psi)) d \mu_\psi \leq  \frac{1}{2}  \int_{\Omega } h(- \psi\slash E_\chi(\psi)) d \mu_\psi  \leq  \frac{1}{2}\cdot
$$
Therefore $\lambda \geq \max \{E_\chi(\varphi) , 2 E_\chi(\psi) c \slash \varepsilon\} \Longrightarrow \int_\Omega h (-\varphi \slash a \lambda) d \mu_\psi \leq 1$ and  we infer
$$
\Vert \varphi \Vert_{ L_h(\mu_\psi)} \leq \max \left\{ a E_\chi(\varphi) , 2^{ 1+1\slash n }  E_\chi(\psi)  \frac{a}{a-1} \right\}.
$$

\medskip

We now no longer assume that $\p \in {\mathcal T}(\Omega)$.
The proof of Theorem \ref{thm:energy-approximation} provides a decreasing sequence $({\psi}_j)$ of functions 
in $\mathcal T(\Omega)$ converging to $\psi$ such that
 $$
 \lim_{j \to + \infty} E_\chi(\psi_j) = E_\chi(\psi)
 \; \; \text{ and } \; \; 
(dd^c \psi_j)^n = {\bf 1}_{A_j} (dd^c \psi)^n,
$$
 where $(A_j)$ is an increasing sequence of subsets converging to $\Omega$.
Applying the previous estimate and taking the limit, we obtain the required inequality 
as we claim that
$ \Vert \varphi \Vert_{ L_h(\mu_\psi)} = \lim_{j \to + \infty} \Vert \varphi \Vert_{ L_h(\mu_{\psi_j})}$.

Indeed  since $j \mapsto \mu_{\psi_j} = (dd^c \psi_j)^n$ is non decreasing and  $\mu_{\psi_j} \leq \mu_{\psi}$,  
the sequence $\Vert \varphi \Vert_{ L_h(\mu_{\psi_j})}$ is non decreasing and 
$\lim_{j \to + \infty} \Vert \varphi \Vert_{ L_h(\mu_{\psi_j})} \leq \Vert \varphi \Vert_{ L_h(\mu_\psi)}$.
On the other hand if $\lambda =   \sup_j \Vert \varphi \Vert_{ L_h(\mu_{\psi_j})}$ then
$\int_\Omega h(-\varphi \slash \lambda) (dd^c {\psi_j})^n  \leq 1$.
Since  
$$
\int_{A_j} h(-\varphi \slash \lambda) (dd^c {\psi})^n = \int_\Omega h(-\varphi \slash \lambda) (dd^c {\psi_j})^n   \leq 1,
$$
the monotone convergence theorem yields $\int_{\Omega} h(-\varphi \slash \lambda) (dd^c {\psi})^n \leq 1$.
We infer $\Vert \varphi \Vert_{ L_h(\mu_{\psi})} \leq  \lambda = \lim_{j \to + \infty} \Vert \varphi \Vert_{ L_h(\mu_{\psi_j})}$,
 proving our claim.
\end{proof}

We have a partial converse to this quantitative integrability property.

\begin{theorem} \label{thm:quantitative}
 Fix $\chi \in {\mathcal W}^+$, $0<a<1$ and a positive Radon measure $\mu$ that satisfies  $I(a,\chi)$.
Then there exists a unique $\psi \in \mathcal E_\chi(\Omega)$ such that $\mu = MA(\psi)$.
\end{theorem}

\begin{proof} 
It follows from \cite[Theorem 4.5]{BGZ09}  that the measure $\mu$ does not charge pluripolar sets.

Assume first that $\mu$ has a compact support $K \subset \Omega$. 
Using a projection argument  \cite[Theorem 6.3]{Ceg98} (see also the proof of 
\cite[Corollary 11.10]{GZbook}), 
we can decompose  $\mu=f MA(u)$ with $u \in {\mathcal T}(\Omega)$ and $0 \leq f \in L^1_{loc}(MA(u))$.  
The bounded subsolution property \cite[Theorem A]{Kol96} ensures that $\mu_j=MA(u_j)$ for some  bounded psh function $u_j$. 
Modifying $u_j$ is a neighborhood of the boundary of $\Omega$ we can further assume that $u_j \in {\mathcal T}(\Omega)$.
Since $j \mapsto \mu_j$ is increasing, 
it follows from Theorem \ref{thm:comparisonprinciple} that $j \mapsto u_j$ is decreasing.
We let $\psi =\lim \searrow u_j$ denote the limit of the $u_j$'s.
By assumption
$$
||u_j||_{L_{\chi}(MA(u_j)} \leq ||u_j||_{L_{\chi}(\mu)} \leq a ||u_j||_{L_{\chi}(MA(u_j)}+C_a,
$$
with $0<a<1$, hence the energies
$$
E_{\chi}(u_j) =||u_j||_{L_{\chi}(MA(u_j)} \leq \frac{C_a}{1-a}
$$
are uniformly bounded. Thus $\psi \in \mathcal E_\chi(\Omega)$
and the continuity of the complex Monge-Amp\`ere operator along decreasing sequences yields
$$
MA(\p)=\lim MA(u_j)=f MA(u)=\mu.
$$
The uniqueness is a consequence of the uniqueness of solutions to such equations in the larger class
${\mathcal E}^1(\Omega)$ (see \cite[Theorem 6.2]{Ceg98}).

\smallskip

When the support of $\mu$ is not compact, we consider an exhaustive sequence of compact sets $(K_j)$ of $\Omega$ and 
consider  $\mu_j= {\bf 1}_{K_j} \cdot \mu \leq \mu$.
By the previous case there exists $\varphi_j \in \mathcal E_\chi(\Omega)$ such that 
$\mu_j = (dd^c \varphi_j)^n$ on $\Omega$.
We show as above shows that the sequence $(\varphi_j)$ decreases to a function $\varphi\in \mathcal E_\chi(\Omega)$ 
such that $\mu = (dd^c \varphi)^n$.
\end{proof}

\subsection{A conjectural integral characterization}

\subsubsection{Rescaling}

Observe that for all $f \in L_{\chi}(\mu)$, one has
$$
||f||_{L_{\chi}({\mu}/{2})} \leq  ||f||_{L_{\chi}(\mu)},
$$
with strict inequality unless $f=0$. It is thus natural to consider
$$
\kappa(\mu,\chi):=\sup_{\f \in {\mathcal E}_{\chi}(\Omega) \setminus \{0\}} \frac{||\f||_{L_{\chi}(\mu/2)}}{ ||\f||_{L_{\chi}(\mu)}},
$$
when ${\mathcal E}_{\chi}(\Omega) \subset L_{\chi}(\mu)$.
The quotient is invariant under rescaling $\f \mapsto \gamma \f$, hence we can restrict
to a relatively compact  subset  of functions
whose $\chi$-energy remains uniformly bounded  and away from zero
(see Proposition \ref{pro:energycpct}).

  \begin{conj} \label{conj:kappa}
We conjecture that $\kappa(\mu,\chi)<1$ for any positive Radon measure such that ${\mathcal E}_{\chi}(\Omega) \subset L_{\chi}(\mu)$.
 \end{conj}
 
 When $\chi(t)=-(-t)^p/p$ it is easy to check that $\kappa(\mu,\chi)<2^{-1/p}<1$, but estimating 
  $\kappa(\mu,\chi)$ for weights with (super)exponential growth seems more delicate.

 \subsubsection{Bedford's problem}
 
 Conjecture \ref{conj:kappa} has  interesting consequences.
 
 \begin{thm} \label{thm:range}
 Assume Conjecture \ref{conj:kappa} holds and fix $\mu$ a positive Radon measure.
 \begin{enumerate}
 \item Fix $\chi \in {\mathcal W}^+$. Then $\mu=MA(\p)$ for $\p \in {\mathcal E}_{\chi}(\Omega)$ if and only if 
 $\mu$ satisfies $I(a,\chi)$ for some  $a \in (0,1)$.
 
 \item $\mu=MA(\p)$ for $\p \in \tilde{{\mathcal T}}(\Omega)$ if and only if 
 for any   weight $\chi \in {\mathcal W}^+$, there exists $0<a<1$ such that 
 $\mu$ satisfies $I(a,\chi)$.
 \end{enumerate}
  \end{thm}
  
 We set here $\tilde{{\mathcal T}}(\Omega)={\mathcal E}^1(\Omega) \cap L^{\infty}(\Omega)$.
 This set consists of bounded plurisubharmonic functions that have zero boundary values
 in the  ${\mathcal E}^1(\Omega)$-sense. This is slightly weaker than
 the condition imposed for test functions, hence ${\mathcal T}(\Omega) \subset \tilde{{\mathcal T}}(\Omega)$.

  \begin{proof}
If  $\mu$ satisfies $I(a,\chi)$ for some $0<a<1$,
it follows from Theorem \ref{thm:quantitative} that  $\mu=MA(\p)$ for 
some $\p \in {\mathcal E}_{\chi}(\Omega)$. 
Assume conversely that  $\mu=MA(\p)$ for some $\p \in {\mathcal E}_{\chi}(\Omega)$.
Thus $2\mu=MA(\p_2)$ with $\p_2=2^{1/n} \p \in {\mathcal E}_{\chi}(\Omega)$.
Using Proposition \ref{pro:quantitative} we can find $1<a< 1/\kappa(2\mu)$ and $C_a>0$ such that  
for all $\f \in {\mathcal E}_{\chi}(\Omega)$,
$$
||\f||_{L_{\chi}(2\mu)} \leq a E_{\chi}(\f)+C_a.
$$
We infer
$$
||\f||_{L_{\chi}(\mu)} \leq \kappa(2\mu) ||\f||_{L_{\chi}(2\mu)} \leq
a \kappa(2\mu) E_{\chi}(\f)+ \kappa(2\mu)C_a.
$$
Thus $\mu$ satisfies $I(a',\chi)$ for $a'=a \kappa(2\mu,\chi)<1$, proving (1).

\medskip

We now prove (2). Assume $\mu=MA(\p)$ for $\p \in \tilde{{\mathcal T}}(\Omega)$.
Since 
$$
 \bigcap_{\chi \in {\mathcal W}^+} {\mathcal E}_{\chi}(\Omega) \subset \tilde{{\mathcal T}}(\Omega),
$$
 it follows from (1) that  for any   weight $\chi \in {\mathcal W}^+$, there exists $0<a<1$ such that 
 $\mu$ satisfies $I(a,\chi)$.
Assume conversely that this property is satisfied for all $\chi \in {\mathcal W}^+$.
Thus $\mu=MA(\p)$ for $\p \in {\mathcal E}_{\chi}(\Omega)$.
By uniqueness of the solution to the Dirichlet problem in ${\mathcal E}^1(\Omega)$, we obtain
$\p \in  \cap_{\chi \in {\mathcal W}^+} {\mathcal E}_{\chi}(\Omega)$, and the conclusion will follow if
we can show that $\p$ is actually bounded.

We prove this by contradiction. If $\p$ is not bounded then the sublevel sets $(\p<-t)$ are non empty for all $t>0$,
and Lemma \ref{lem:classic} (or rather its extension to ${\mathcal E}^1(\Omega)$, see Remark \ref{rem:keylemE1})
 ensures that
$$
\e: t \in \R^+ \mapsto \mu(\p<-t) \in \R^+_*
$$
is (strictly) positive, and decreases to zero at infinity. We let $h$ denote the Young weight such that
$h(0)=0$ and $h'$ is defined by
$$
h'(t)=\frac{1}{\e(2^j t)}
\; \; \text{ if } \; \; j \leq t < j+1.
$$
By construction we obtain $\int h(\p/\la) MA(\p)=+\infty$ for all $\la >0$. 
Indeed 
$$
\int h(\p/\la) MA(\p)=\int_0^{+\infty} h'(t/\la) \mu(\p<-t) dt
\geq \int_j^{+\infty} \frac{1}{\e(2^j t / \la)} \e(t) dt =+\infty
$$
if $2^j \geq \la$, since $\frac{1}{\e(2^j t / \la)} \e(t) \geq 1$.
Thus $\p \notin {\mathcal E}_{\chi}(\Omega)$.
\end{proof}

  The second point   provides a satisfactory answer to Bedford's Problem \ref{pbm:bedford},
  contingent upon the validity of Conjecture \ref{conj:kappa}.
  This problem remains open since the foundational works of Bedford-Taylor,
  more than forty years ago.

\end{document}